\documentclass{article}

\usepackage{arxiv}

\usepackage[utf8]{inputenc} 
\usepackage[T1]{fontenc}    
\usepackage{hyperref}       
\usepackage{url}            
\usepackage{booktabs}       
\usepackage{amsfonts}       
\usepackage{nicefrac}       
\usepackage{microtype}      
\usepackage{lipsum}
\usepackage{graphicx}
\graphicspath{ {./images/} }
\usepackage{amssymb}

\usepackage{graphicx}
\usepackage{color}
\usepackage{epstopdf}
\usepackage[english]{babel}
\usepackage[export]{adjustbox}
\usepackage{amsmath, amssymb, amsfonts}

\title{Optimal vaccine roll-out strategies with respect to social distancing measures for SARS-CoV-2 pandemic}

\author{
 Konstantinos Spiliotis\\
  Institute of Mathematics\\ 
  University of Rostock\\
D-18057 Rostock, Germany\\
  \texttt{konstantinos.spiliotis@uni-rostock.de} \\
   \And
    Constantinos Ch. Koutsoumaris \\
  Department of Civil Infrastructure \\
  and Environmental Engineering\\
  Khalifa University of Science and Technolog\\
   P.O. Box 127788\\
Abu Dhabi, United Arab Emirates\\
  \texttt{konstantinos.koutsoumaris@ku.ac.ae} \\
  \And
 Andreas Reppas \\
Charité-Universitätsmedizin Berlin\\
corporate member of Freie Universität Berlin\\
Humboldt-Universität zu Berlin\\
Berlin, Germany\\
  \texttt{andreas.reppas@charite.de} \\
  \And
   Jens Starke\\
  Institute of Mathematics\\ 
  University of Rostock\\
D-18057 Rostock, Germany\\
  \texttt{jens.starke@uni-rostock.de} \\
   \And
    Haralampos Hatzikirou $^1$ \\
  Mathematics Department\\ 
  Khalifa University of Science\\
  and Technology, P.O. Box 127788\\
Abu Dhabi, United Arab Emirates\\

}

\begin{document}
\footnote{Corresponding author} 
\maketitle
\begin{abstract}
Non-pharmacological interventions (NPIs), and in particular social distancing, in conjunction with the advent of effective vaccines at the end of 2020, aspired for the development of a protective immunity shield against the spread of SARS-CoV-2. The main question rose is related to the deployment strategy of the two doses with respect to the imposed NPIs and population age. In this study, an extended (SEIR) agent-based model on small- world networks was employed to identify  the optimal policies against Covid 19 pandemic, including social distancing measures and mass vaccination. To achieve this, a new methodology is proposed to solve the inverse problem of calibrating an agent’s infection rate with respect to vaccination efficacy. The results show that deploying the first vaccine dose across the whole population is sufficient to control the epidemic, with respect to deaths, even for low number of social contacts. Moreover, for the same range of social contacts, we found that there is an optimum ratio of vaccinating ages $>65$ over the younger ones of 4/5.
\end{abstract}

\section{Introduction}

For more than a year after the first case of coronavirus disease-2019 (COVID-19),the world is still experiencing the largest social isolation in history \cite{Chou20}. Namely, the Non-pharmacological interventions (NPIs) have been considered as the main strategy for hindering the spread of severe acute respiratory syndrome coronavirus 2 (SARS-CoV-2). Social-distancing, stay-at-home and lockdowns are examples of such interventions targeting at reducing health system overload and overall death toll \cite{Ben21,Haug20}. Nevertheless, NPI consequences on social and economic life raise questions on the impact of these restrictions to the society \cite{Ben21, Haug20}. 

On the other hand, vaccine introduction against SARS-CoV-2 set a new roadmap on fighting the virus spread. By the end of 2020, a variety of vaccines having different efficacy on establishing individual immunity against the virus, were introduced \cite{maier21}. Thus, a vaccination process could accelerate the way towards building an effective safety wall against the virus transmission.   

The vaccination process varied widely among countries. Israel, which manage to reach a vaccination daily rate of over 1\% of its population, became the first country to achieve a complete, doubled-dose vaccination on 65\% of the adult population by April of 2021 using Pfizer/BioNTEch's vaccine \cite{web:Israel}. Such vaccination process resulted in very low numbers of daily new cases and hospitalized patients while the mortality rate due to COVID-19 aftermath was decreased \cite{web:Israel}. Following a similar policy, Great Britain provided a single vaccine dose of AstraZeneca vaccine, with a preference towards older people, and achieved a daily rate close to 1\% \cite{web:UK}. This vaccination process resulted in blocking the epidemic wave started in autumn of 2020 fading out with low numbers of daily new cases by April of 2021. On the other hand, countries like Germany and Greece, that achieved daily vaccination rates of around 0.5\% the period before April 2021, confronted a new exponentially increasing wave of daily new cases started on January-February of 2021 \cite{web:Germ,web:Gr}.

\begin{figure}
    \centering
    \includegraphics[width=0.7\textwidth]{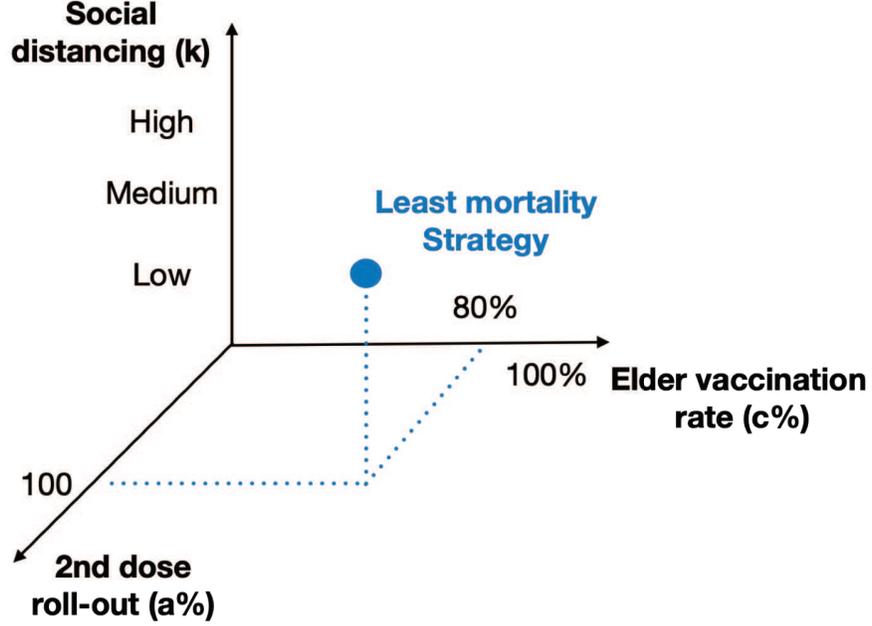}
    \caption{Combining policies: the single / double vaccine dose, the age stratification and the social distancing, in order to mitigate the  SARS-CoV-2 impact. The optimal combination (blue marker) of policies leads to a minimum number of expected deaths.}
    \label{fig:fig1Hatz}
\end{figure}

The above observations show the importance of implementing a carefully designed vaccination protocol over the entire population. However, constraints on the vaccine availability as well as the limited capability of health systems to offer adequate vaccination rates rise the question of how a country can design a successful vaccination plan. For the latter restrictions is under consideration which policy is the best: 1) less people with high resistance to virus (double doses) or 2) more vaccinated people with smaller resistance to virus (singe dose for all). It is an open problem the suitable combination of single or double dose strategy that results in the lower expected number of deaths (see also Fig. \ref{fig:fig1Hatz}).

 At the same period, the vaccination age question arises: The vaccine administration should be from gross elder to younger person, or to start from the beginning with a combination between different group ages. Obviously, older vaccinated people decrease their mortality. One the other hand, younger individuals with increased mobility and less precautions are potential powerful transmitters of virus. Decreasing the infection in this group should have strong impact in the epidemic dynamics and consequence in the expected number of deaths.

Recent studies support that a vaccination planning cannot be successful in the absence of NPIs \cite{Moo21}. In a similar pace, serious concerns has been raised on vaccine efficacy on old individuals \cite{Matr21} or because of the virus mutations \cite{Wis21}. Moreover, the complications followed SARS-CoV-2 virus infection especially to old people suggests that a vaccination policy should also take into account the age distribution of each population \cite{Log21,Log20}, see also Fig. \ref{fig:fig1Hatz}

Mathematical modelling provides a solid ground on studying virus dynamics like the ones observed in the case of SARS-CoV-2 and to extract general disease characteristics, like the speed of the virus spread and the long-term epidemic behaviour \cite{Siet15,Rep10,Foy21,Matr21,Syg20,maier21,Anas20}. Mathematical- computational approaches model the epidemic dynamics for separated groups (e.g. healthy susceptible, infected and recovered) using two approaches: a continuous formulation in the form of differential equations (DE) or an agent based model, where epidemic evolution results from individuals properties and rules \cite{Fer05,Fer06,Siet15,Zha15,Syg20,Gior20,Matr21,maier21,Anas20}: Similar to DE's, agent based models describe the macroscopic population-scale evolution, as this resulting from detailed knowledge on the  microscopic level e.g. the immune mechanisms, host-microbe and host-host interactions, incubation period and mortality \cite{Rep10}. 

The use of networks in the epidemiology is obviously justified by viewing the infection transmission as social contact problem \cite{Eam15,Zha15,Sal10,Rep15}. The advantages of network relay on many topics with most important the heterogeneity of connectivity structure: 
Many epidemic models (e.g. the DE modelling) assumes the well-mixed population, but in real life this is not the case. Each person has distinct number of contacts whom maybe infects \cite{Keel05,Keel05b,Eam15}. For this reason the knowledge of network heterogeneity and its impact on epidemic dynamics are fundamental for predicting and designing NPI's with or without vaccination process for controlling the outbreak of diseases. Poor understanding of the evolving macroscopic dynamics due to contact interactions through the epidemic network may have to serious and fatal consequences \cite{Ioa21}.

In this study, we evaluate the effect of different vaccination policies in the context of different vaccination dose strategies, social distancing and vaccination targets based on the age distributions, see Fig. \ref{fig:fig1Hatz}. By using a network-based epidemic model, we show the interdependence of the above main factors in designing an effective vaccination policy. Initially, the epidemic model is presented. Then, a rigorous methodology to calibrate individual infection rates with respect the vaccination efficacy is proposed. This is achieved be using an equation free methodology by coupling agent based model with numerical solver scheme like Newton-Raphson  \cite{Kev09,Siet15,Rep10}. We notice that the multiscale inverse problems are particularly challenging and remain on active research field \cite{Chae2020,GOE15,MARIN2020}. Next, we show how the network structure, which serves as a proxy of the NPIs, affects the epidemic output by designing different vaccination plannings. Finally, we analyze the vaccination efficacy by considering both the age distribution for selecting the vaccination targets along with the use of NPIs.

The results show that a single vaccination dose is sufficient to control the epidemic when mild social distancing (medium number of contacts) are implemented. An optimum behaviour exists for age vaccination planning combined with the implementation of mild NPIs. Therefore, we can state that the implementation of mild NPIs support any vaccination planning towards controlling the epidemic effect and building the desired herd immunity wall against the virus spread.

\begin{figure}[t!]
\begin{center}
\begin{picture}(400,300)
\includegraphics[width=13cm]{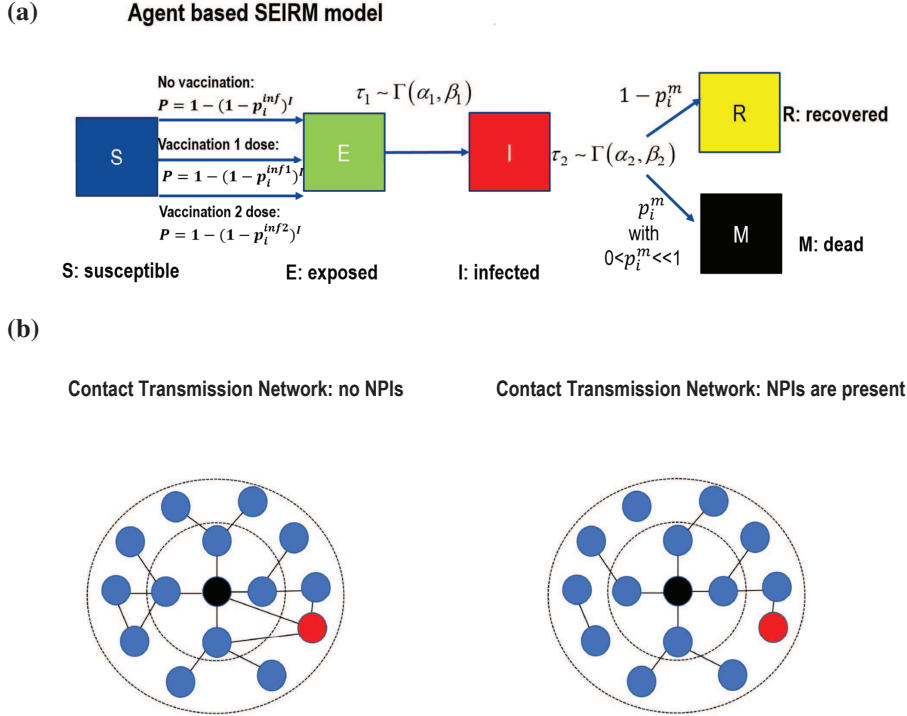}
\put(-365,270){\textbf{(a)}}
\put(-365,150){\textbf{(b)}}
\end{picture}
\end{center}
\caption{\textbf{(a)} Representative figure of the model. A susceptible agent (S) becomes  exposed (E) with probability $P$, see eq. \eqref{eq:inf_prob}. The  probability $P$ depend on the infection rate $p_{i}^{inf}$. Vaccinated agent shows higher resistance to infection by decreasing the infection rate : $p_{i}^{inf2}<p_{i}^{inf1}<p_{i}^{inf}$ where $p_{i}^{inf1}$ and $p_{i}^{inf2}$ stands for single/ double dose respectively.Then, after time $\tau$ chosen from $\Gamma$ distribution the exposed become infected (I). Finally the agent is recovered (R) with probability $(1-p_{i}^m)$, $0<p_{i}^m<<1$ or pass away(M). \textbf{(b)} For representative reasons we show a simplification of contact networks in the absence  and second, in the existence  of NPI's. In the later case, main difference is the reduction of contacts which corresponds on the reduction of k-degree  parameter in the small world structure.  
}
\label{fig:Fig1}            
\end{figure}

\section{Materials and Methods}
\subsection{Agent based epidemic model on network for designing vaccination policies on NPIs}
For studying the SARS-CoV-2 epidemic dynamics a discrete agent-based on complex network is used. The model is a SEIRM as introduced in \cite{Anas20,Foy21,Cal20}. Each node represents a distinct agent with indirect edge as social contacts. The network is constructed according to the following simple way by: Initially, connecting every node with its $k$ nearest neighbors in a ring-like topology and secondly, rewiring every link to a random node with probability $p$ \cite{Wat98}. Each node-individual can have 5 states S,E,I,R,M susceptible or healthy, exposed to virus, infected asymptomatic or not, recovered and cannot be infected any more and unlucky dead (mortality). Every time step corresponds to a day. Depending on the vaccination roll-out policy, a specific constant number of susceptible individuals (per day),randomly chosen administrate the first or second dose, enhanced their resistance to SARS-CoV-2 infection. Each $i-$th healthy -susceptible individual  becomes exposed with probability 
\begin{equation}
    P_i = 1-(1-p_{i}^{inf})^{I_i}
    \label{eq:inf_prob}
\end{equation}
where $I_i$ is the number of infected neighbours of the $i$ susceptible. The $P_i$ is an increased function with respect the variable $p_{i}^{inf}$, with $P_i(0)=0, P_i(1)=1$. Non-vaccinated susceptible agents have the same constant value $p_{i}^{inf}$, while vaccinated individuals decrease the $p_{i}^{inf}$ to $p_{i}^{inf1}$ for the first dose (after 10 days) and  to $p_{i}^{inf2}$ if second dose applied (after 30 days, from first dose), with $p_{i}^{inf2}<p_{i}^{inf1}<p_{i}^{inf}$, the new values of $p_{i}^{inf}$ depends on the efficacy of the vaccine. In sec. \ref{sec:Newton} a rigorous methodology to extract these values as function of vaccine efficacy is provided. Upon infection, the agent changes state from susceptible to exposed.  Then a waiting time is assigned $\tau_E$ which is drawn from $\Gamma$ distribution i.e. $ \tau_E \thicksim \Gamma(a_E,b_E)$ where $a_E, b_E$ denote the parameters of $\Gamma$ distribution.
During every time step, the waiting time $\tau_E$ 
is reduced by 1, i.e. $\tau_E=\tau_E-1$ and $\tau_E>0$. The first time step where  $\tau_E<0$ the disease progresses, and the exposed individual becomes infected. Then a new waiting time is assigned,  $\tau_I \thicksim \Gamma(a_I,b_I)$, from $\Gamma$ distribution. At its time step, again the $\tau_I$  is decreased by 1, i.e.  $\tau_I=\tau_I-1$. In the model, vaccinated people can still  be infected which implies that they can also still infect other individuals, but with lower transmission rate.

Finally, the transition from infection to recover or unlike to death is depended on the age and on administration of vaccination. Specifically, in the model the first time when $\tau_I<0$ the node changes state according to the following way: With small probability $p_i^m$ the infected dies (mortality) or in the opposite case (with probability $1- p_i^m$) the agent recovers and cannot be exposed or infected any more. The  $p_i^m$ value is defined according to the following way:
If the $i$-th individual is less than 65 years old then independently of the vaccination or has probability is  $p_i^m=0.005$. In the opposite case i.e. age $>$ 65 then the non-vaccinated people have probability $p_i^m=0.13$, while in the vaccinated case the $p_i^m$ decreases to $p_i^m=0.008$.

A sketch for the SEIRM dynamics can be found in Fig.~\ref{fig:Fig1}(a).
The parameters of $\Gamma$-distributed waiting times inferred from COVID-19 disease characteristics, as these have been found to describe the disease best \cite{Lin20,Syg20}.  Specifically for the exposed state $\tau_E \thicksim \Gamma(a_1,b_1)$ and $\tau_I \thicksim \Gamma(a_2,b_2)$ ,where $a_1=9, b_1=1/3$ and  $a_1=100/3, b_1=3/10$. 
For these parameters of $\Gamma$ the mean value and the standard deviation is 3 and 1 and 10 and $\sqrt{3}$, respectively. 
In the next section a rigorous methodology is proposed in order to estimate the new infection probability $p_{i}^{inf}$  value for vaccinated individuals.



\subsection{Simulating the infection resistance with respect vaccination efficacy}
Vaccinated individuals increase their resistance to infection, which is implies lower value of $p_{i}^{inf}$.  The new values of  $p_{i}^{inf}$  i.e. $p_{i}^{inf1}$ or $p_{i}^{inf2}$ will depend on the efficacy of administrated vaccine and it constitutes an open problem for agent based epidemic model the definition of consistent value of $p_{i}^{inf}$. The idea behind our approximation for estimating the new $p_{i}^{inf}$ in consistent with vaccine efficacy is that we can mimic the physical process (see sec.~\ref{sec:clinic}) of efficacy estimation \cite{Log21} using the agent based model as timestepper or map \cite{Kev09,Mars14,Rep10,Spil10}, with input the unknown $x=p_{i}^{inf}$ and output $d=\Phi_T(x)$, the reproduced efficacy. Then using the equation free methodology \cite{Kev09,Mars14,Rep10,Spil10,Mars14} an equation solver (i.e. Newton Rapshon scheme) is  coupled with the map  to find  numerically the demanded value  $x=p_{i}^{inf}$.
\subsubsection{Vaccine efficacy in clinical trials}
\label{sec:clinic}
In clinical trials e.g. \cite{Log20,Log21} in order to define a vaccine's efficacy, the following procedure is applied: Individuals from a human sample, administrated with true vaccination or placebo. Then after a time period, the odds ratio (OR) is defined as the  fraction of the proportion of vaccinated people who developed COVID-19 (from the true vaccinated ensemble), divided  by the proportion of COVID-19 developed from the sample of placebo. Specifically, if $a$ is the  number of vaccinated participants which COVID-19 developed and $b$ is the number of vaccinated participants. The ratio $\frac{a}{b}$ is the proportion of COVID-19 in the sample of vaccinated people. Clearly higher vaccine efficacy  resulting to low values of the ratio $\frac{a}{b}$ i.e. $\frac{a}{b}\rightarrow0$. Similar if $c$ is the number of placebo participants with COVID-19, and $d$ is the total number of placebo then the  $\frac{c}{d}$ is the percentage of individual with COVID-19 in the sample of placebo. The odds ratio (OR) is given by
\begin{equation}
  OR=\frac{{a}/{b}}{{c}/{d}}
\end{equation}
then, the efficacy $e$ of the vaccination or the efficacy is defined as
\begin{equation}
    e = 1-OR
\end{equation}
If the vaccine is highly effective then $\frac{a}{b}\rightarrow 0 \implies \ OR\rightarrow 0 \implies e \rightarrow 1-0=1$. Instead, for low efficacy implies $\frac{a}{b}\approx \frac{c}{d} \implies OR\approx 1 \implies  e = 1-OR=1-1=0$.
For instance, Sputnik vaccine first dose (rAd26) and the second dose (rAd5) \cite{Log21}, from almost 21000 participant finally 14964 and 4902 received the second dose of vaccine and placebo respectively and were included in primary outcome analysis \cite{Log21}. According to  \cite{Log21}, 21 days after the 1 dose (tab.2 in there in) 16 and 62 respectively from the vaccine and placebo sample were infected. Then the $OR=\frac{{14}/{14964}}{{62}/{4902}}=0.0845$ resulting to $e=1-0.0845=0.915$.


\begin{figure}[t!]
\begin{center}
\begin{picture}(400,250)
\includegraphics[width=13cm]{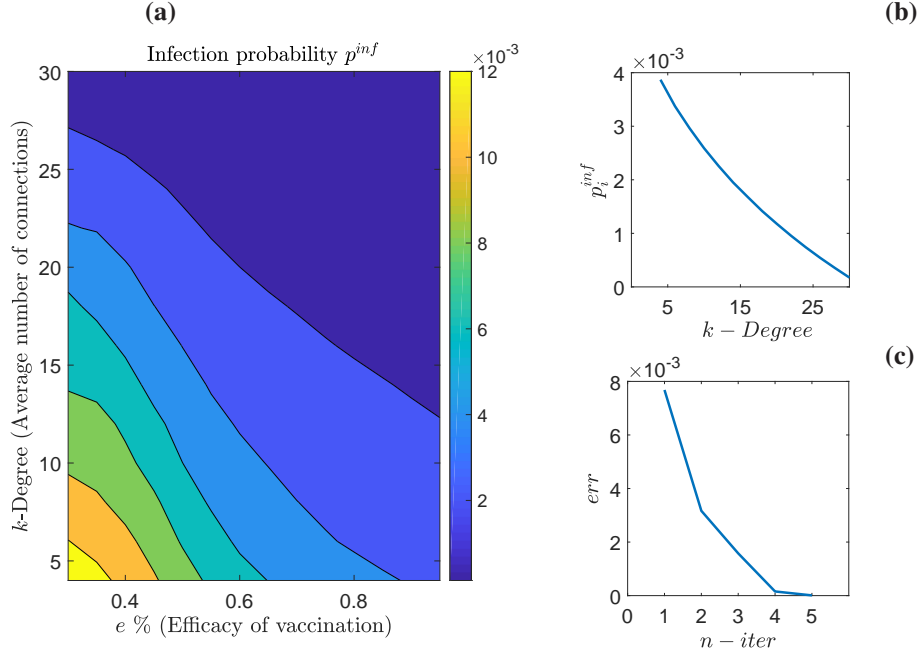}
\put(-300,245){\textbf{(a)}}
\put(-20,245){\textbf{(b)}}
\put(-20,115){\textbf{(c)}}
\end{picture}
\end{center}
\caption{\textbf{(a)} The values of $p_{inf}$ with respect the efficacy $e_0$ and the mean $k$-degree (number of social contacts), resulting from eq.\ \eqref{eq:Equa_Eff}. \textbf{(b)} The values of $p_{inf}$ for $e_0=0.9$ resulting from eq. \ \eqref{eq:Equa_Eff} for all degrees. \textbf{(c)} The error of Newton Raphson using eq. \eqref{eq:deriv} for $k=10$.}
\label{fig:Newton}            
\end{figure}

\subsection{Inverse problem methodology\text{:} estimating the individual infection probability for a given vaccine efficacy}
\label{sec:Newton}
Here we address the question, given a vaccination efficacy how can we determine the individual infection probability. Depending on the vaccine efficacy we expect a decrease in the probability of infection or more precisely the transition from $S\rightarrow E$, i.e. the $p_{i}^{inf}$ in eq.\ \eqref{eq:inf_prob}.

In the computational agent based model, two samples for vaccination and placebo  are randomly chosen. Each vaccinated agent alters the initial infection value $p_{i}^{inf}$ to a new unknown value $p^{inf}$ (i.e. $p_i^{inf1}$ or $p_i^{inf2}$  ) where this value has to be determined.  At the end of short time period of 100 step(days) the ratio $OR$ and the efficacy $e = 1-OR$ are defined in the exact same way as in sec. \ref{sec:clinic}. The final value resulting from average on the ensemble $N_{times}$ identical repetitions.

For a constant structure-network, the above procedure defines a map (or timestepper) $\Phi_T:[0, 1]\rightarrow \mathbb{R}$  \cite{Rep10,Rep15,Kev09} :
\begin{equation}
    e = \Phi_T(A,p^{inf})
    \label{eq:d_phi}
\end{equation}
where $A$ is the adjacency matrix of the network and contains all the information about the structural topology (degree distribution, clustering, path length etc) and the variable $p_{i}^{inf}$ has to be determined  in order to achieve efficacy $e = e_0$. Then, the following equation has to be solved 
\begin{equation}
    e_0=\Phi_T(A,p^{inf})\Leftrightarrow e_0-\Phi_T(A,p^{inf})=0 \Leftrightarrow G(p^{inf})=0 
    \label{eq:Equa_Eff}
\end{equation}
The above equation can be solved numerically for example using the  Newton Rapshon algorithm.
The derivative of $G$ can be approximated using a numerical approximation scheme i.e.
\begin{equation}
    G'(p^{inf})\approx \frac{G(p^{inf}+dp)-G(x)}{dp}= \frac{\Phi_T(A,p^{inf}+dp)-\Phi_T(A,p^{inf})}{dp}
    \label{eq:deriv}
\end{equation}
The solution $p^{inf}$ for different values of $e_0$ defines the inverse function of eq.\ \eqref{eq:d_phi} and resulting in the value $p_{i}^{inf1}, p_{i}^{inf2}$.

Of course these values depend on the network structure (i.e. degree $k$ and rewiring probability $p$). Fig.~\ref{fig:Newton}(a) shows the values of $x=p_{i}^{inf}$ with respect the different levels of efficacy $e$ and connectivity- degree $k$, (with constant rewiring probability $p=0.005$). Fig.~\ref{fig:Newton}(b) shows the decreasing monotonously of $x=p_{i}^{inf}$ with respect the local degree $k$. Fig.~\ref{fig:Newton}(c) depicts the convergence of Newton-Rapshon scheme for eq.\ \eqref{eq:Equa_Eff}, using the approximation derivative i.e. eq. \eqref{eq:deriv} and for $e_0=0.9, k=10$.

\begin{figure}[t]
\begin{center}
\begin{picture}(400,280)
\includegraphics[width=12cm]{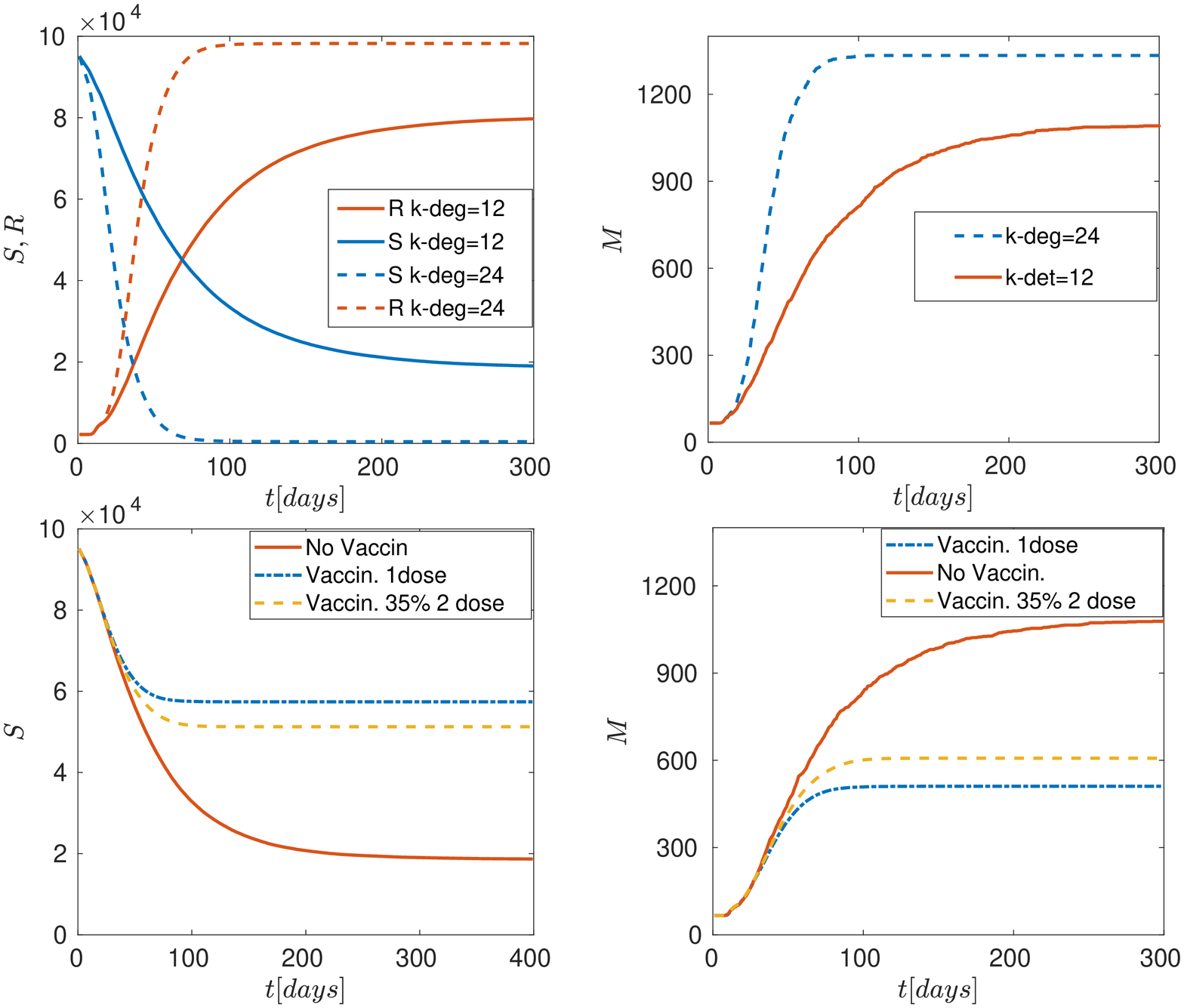}
\put(-340,280){\textbf{(a)}}
\put(-180,280){\textbf{(b)}}
\put(-340,130){\textbf{(c)}}
\put(-180,130){\textbf{(d)}}
\end{picture}
\end{center}
\caption{Epidemic dynamics for different NP policies  \textbf{(a)}Without vaccination. Time series of susceptible (S) and recovered (R) for different local connectivity k=12  and k=24.  \textbf{(b)} The expected number of deaths (without vaccination) for different local connectivity, $k=12$  and $ k=24$. \textbf{(c)} and \textbf{(d)} Epidemic dynamics for different vaccination policies. Solid line corresponds to no vaccination, dash-dot line to one dose and dash line to $a=35\%$ under the assumption of constant number of vaccination w=1000 per day.\textbf{(c)}  Time series of susceptible (S)  \textbf{(d)} The expected number of deaths for local connectivity k=12, reveals the more efficient policy of one dose.}
\label{fig:existence}            
\end{figure}


\section{Results}
For all simulations a total number of $N=10^5$ agent (the population) was used. The probability of rewiring probability  $p=0.005$ was selected in the small world construction. Different dosing rates: 10, 500 and 1000 doses/day, correspond to a percentage $0.01\%, 0.5\%$ to $1\%$ of our population are studied. The infection probability  $p_{i}^{inf}$ in case of vaccination is adapted to each efficacy and local connectivity $k$ according to theory of sec.\ref{sec:Newton}. An unvaccinated agent/person has  $p_{i}^{inf}=0.02$ \cite{Syg20}, which also results to {$e=0$}. 

Initially, the single or double dosing strategy investigated, assuming a constant number of vaccination per day. 
Since  vaccination involves two doses separated three-four weeks apart, the first question rose is how the authorities should vaccinate the people in order to impact the virus spread and to decreases to a minimum number of deaths, combining the social distancing measures (low-medium-high number of contacts). To deal with this problem a percentage ratio $a\%$ is defined on the daily number of second dose uptake. For example in case of $100$ vaccination per day, a ratio  $a=10\%$ means that 90 new vaccinations will be performed (1st dose) and the remain 10 will take the second dose per day. In the following simulations, single dose results to efficacy of 50\% while with the double dose (full vaccinated) the efficacy is 90\%.
 
The second approach estimates the age priority strategy during single dose roll out policy for medium and high efficacy vaccination process, consider a combination of different social distancing measures (low-medium-high number of contacts). For simplicity, two age populations groups are defined. The group A includes individuals with age over 65 and represent the 20\% of the population and the group B includes all the people less than 65 years old. The ratio $c\%$ is defined to express the priority between the two age groups. Specifically, $c\%$ gives the percentage of the A group greater than 65 years old in which are daily vaccinated. For example in constant window of $100$ daily injections, if $c=0.3$, then 30 individuals belong to group A (greater than 65 years old) and the rest 70  belong to group B.  

 In all simulations the expected number of deaths results as average value of 128 identical realization of the model. Fig.~\ref{fig:existence} depicts four different scenarios for a constant rate of 1000 doses/day (or the 1\% of population). The simplest case of the epidemic dynamics without vaccination is showed in Fig.~\ref{fig:existence}(a)(b). The well-documented problem of fast dynamics ($k=24$ in  Fig.~\ref{fig:existence}(a)(b)) due to big number of contacts is recapitulated. The susceptible population in case of $k=24$ (dash line) goes to zero before 100 days and similarly the expected number of deaths goes to around 1300 Which is considered to be high, making the health care system incompetent to react. In Fig.~\ref{fig:existence}(c),(d) shows the vaccination effects for different policies for $k=12$: single dose (dash dot in Fig.~\ref{fig:existence}(c),(d)),  65\% per day admits the first dose and the remain 35\% the second dose (i.e. in w=100 per day, 65 are new vaccinated individuals and 35 admit the 2nd dose) also the no vaccination is depicted for comparison reasons. The results strongly suggest the investigation optimal strategies that combines social distancing restrictions and vaccination process.

\begin{figure}[t!]
\begin{center}
\begin{picture}(400,250)
\includegraphics[width=13cm]{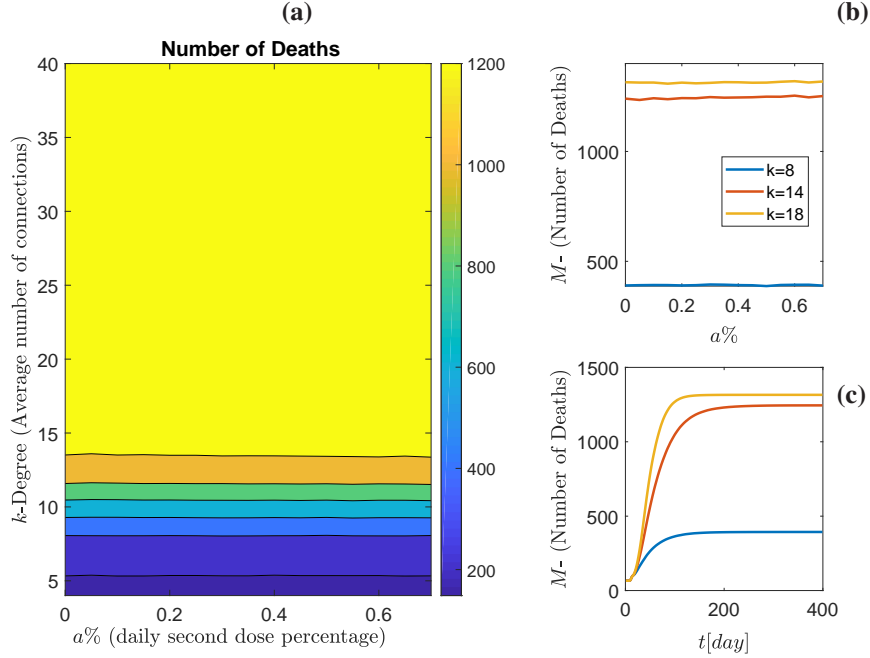}
\put(-240,245){\textbf{(a)}}
\put(-30,245){\textbf{(b)}}
\put(-30,100){\textbf{(c)}}
\end{picture}
\end{center}
\caption{Simulation of epidemic model with respect the administration of dose strategy. Each day a constant ratio of 0.01\% of the population per day (10 people) are vaccinated. \textbf{(a)}The number of dead with respect the degree $k$ (social contacts) and percentage a\% of people take the second dose. There are well defined regions of different behaviors with major contribution of highly mortality yellow area. \textbf{(b)}The efficacy of NPI with respect the single or double dose. Three representative degrees $k=8, 14, 18$ are shown.  \textbf{(c)} Time series of the number of death for constant mixture $a=0.35$ and for degree $k=8, 14, 18$. }
\label{fig:Roll10}            
\end{figure}

\subsection{Single/double dose vaccination policy under low daily roll out  injection rate}

Initially, a low number of daily vaccinations is studied. Each day a constant number of 10 people/day is vaccinated, corresponds to 0.01\% of the whole population (roll out rate). The expected number of deaths with respect to the paramater a\% and the NPIs is depicted in Fig.~\ref{fig:Roll10}. In \ref{fig:Roll10} (a) depicts the number of deaths with respect the degree of connections k and the percentage a\% of people take the second dose. 

Clearly, with very low number of vaccinations per day, only the degree restriction (lock down) can reduced the expected number of deaths. Three representative degrees $k=8, 14, 18$ are shown in Fig.~\ref{fig:Roll10}(b) with respect the percentage of single or double doses and the corresponding average time series are depicted in  Fig.~\ref{fig:Roll10}(c). Almost constant lines in  Fig.~\ref{fig:Roll10}(b) reveal the independence of number of deaths with respect  to dosing. In a nutshell, for low vaccination roll out the dosing policy has no any impact to the virus spread and the deaths of it. 

Finally, due the small number of daily vaccinations the dynamics is dominated from social connections and unfortunately for $k>13$ the expected number of deaths stay unchanged to very a high number. The good scenario of low number of deaths is achieved only in a narrow area of 8-9 degrees (contacts), which is very difficult to be practically achieved. 

\begin{figure}[t]
\begin{center}
\begin{picture}(400,260)
\includegraphics[width=13cm]{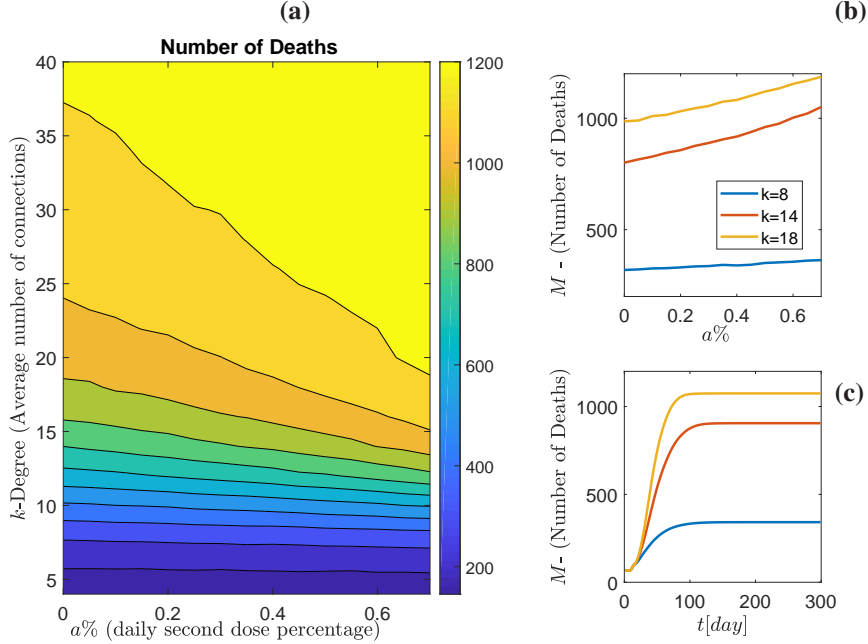}
\put(-240,245){\textbf{(a)}}
\put(-30,245){\textbf{(b)}}
\put(-30,100){\textbf{(c)}}
\end{picture}
\end{center}
\caption{Simulation of Epidemic model with respect the dose strategy. A constant number of 500 people or 0.5\%  of the population (roll out rate) are vaccinated per day. The percentage a\% describes the proportion of people take the second dose and increase the infection resistance. \textbf{(b)}The efficacy of NPI with respect the single or doubled dose. Three representative degrees (social contacts) $k=8, 14, 18$ are shown.\textbf{(c)} Time series of the number of death for constant mixture $a=0.35$ and for degree $k=8, 14, 18$.}
\label{fig:Roll500}            
\end{figure}

\begin{figure}[t]
\begin{center}
\begin{picture}(400,250)
\includegraphics[width=13cm]{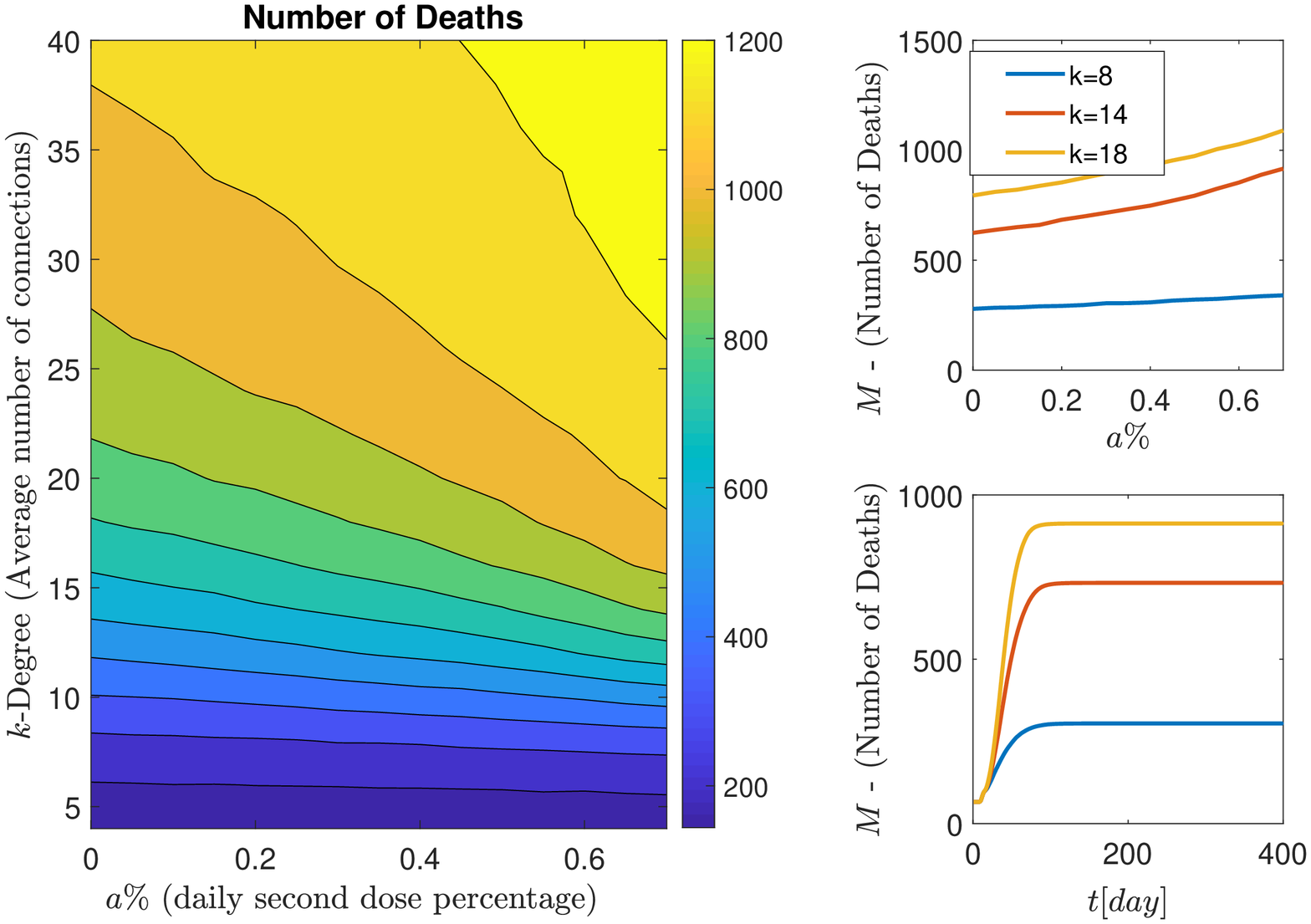}
\put(-240,230){\textbf{(a)}}
\put(-30,230){\textbf{(b)}}
\put(-30,100){\textbf{(c)}}
\end{picture}
\end{center}
\caption{Simulation of Epidemic model with respect the dose strategy. 1000 people or 1\% per day of the population are vaccinated.  \textbf{(b)}The efficacy of NPI with respect the single or doubled dose. Three representative degrees (social contacts)  $k=8, 14, 18$ are shown. \textbf{(c)} Time series of the number of death for constant mixture $a=0.35$ and for degree $k=8, 14, 18$.}
\label{fig:Roll1000}            
\end{figure}

\subsection{A single dose vaccination policy is sufficient to control the epidemic for large roll out rates?}

In this section a different scenario of vaccination policy with increased number of vaccinations is studied. Each day a constant number of 500 people, corresponds to 0.5\% of the population, are vaccinated. The average dynamics (over the example of 128 realization) with respect to the $a\%$ and the connectivity degree are shown in Fig.~\ref{fig:Roll500}(a). In Fig.~\ref{fig:Roll500}(b) the expected number of deaths with respect to a\% for three representative degrees $k=8, 14, 18$ are shown, while in  Fig.~\ref{fig:Roll500}(c) the corresponding time series of the number of death (in average) for constant mixture $a=0.35$ and for the same degree as above $k=8, 14, 18$ are shown. Fig.~ \ref{fig:Roll500}(a) gives an overall description of the behavior by partitioning the domain into zones of expected number of death with respect one /two dose and NPI's policies. A coarse separation of \ref{fig:Roll500}(a), the dynamics results to three zones of low(deep blue), medium(blue-green) and high(yellow) number of deaths. The expected number of deaths still remains monotonically increasing with respect to the $a\%$, meaning that independently of degrees(contacts) the best policy comes with one dose (see also Fig.~\ref{fig:Roll500}(b)). Additionally,  the increase of $a\%$ has high negative impact since it increases dramatically the  expected number of deaths. For example in Fig.~\ref{fig:Roll500}(b) the curve with mean degree of k=14, starts at around 800 (a\%=0) as expected number of deaths and increases monotonically until 1100 for a=0.7. 

Similar behavior is obtained when the daily vaccination rate is 1\%. Again there are three coarse zones of dynamics of low -medium and high number of the expected number of deaths. Of course the higher number of vaccinations decreases the expected number of deaths when compared to the 0.5\% vaccination rate. Comparing the areas of two yellow zones (light normal and dark yellow)in Fig.~\ref{fig:Roll500}(a) and Fig.~\ref{fig:Roll1000}(a) we find that in case of vaccine rate 0.5\% corresponds to 60\% while in case of vaccine rate 1\% corresponds to 40\% of the hall area showing the obvious advantage of increased daily vaccination number. A clear view is the comparison of \ref{fig:Roll500}(c) and Fig.~\ref{fig:Roll1000}(c), where e.g. for $k=18$ the mortality is around 1074 and 912 respectively a reduction of 15\%., while for  for $k=14$ the mortality is around 905 and 732 respectively a reduction of 19\%.

Importantly, the single dose policy has also great social impact in the sense that even with high contact number i.e. degree $15<k<22$ (no restrict lock down, see also \cite{Ioa21}), the expected number of deaths remains low (while for example with $k=20$ and $a=0.6-0.7$ the expected number is $\approx$ 1000, passing in yellow zone). To conclude, the choice of double dose in a constant roll out rate can increase dramatically the expected number of deaths. 

\begin{figure}[ht]
\begin{center}
\begin{picture}(400,250)
\includegraphics[width=13.5cm]{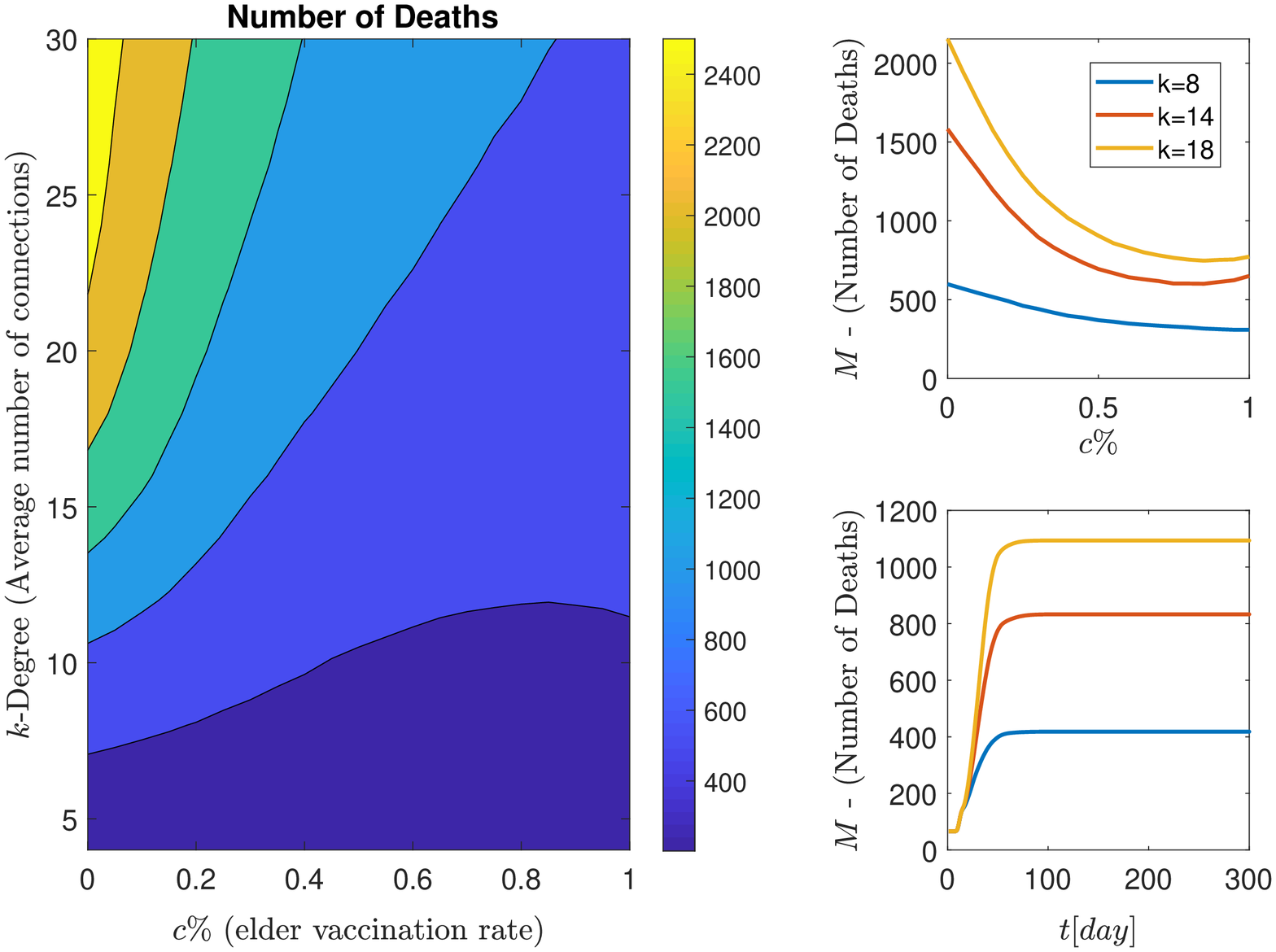}
\put(-240,250){\textbf{(a)}}
\put(-30,250){\textbf{(b)}}
\put(-30,100){\textbf{(c)}}
\end{picture}
\end{center}
\caption{Simulation of epidemic model with respect the age priority and with combination of NPI's. Each day a constant number of 1000 people/day (roll out window) are vaccinated with one dose of $e_0=0.7$ efficacy. The percentage c\% describes  the number of above 65 year which are vaccinated.  \textbf{(a)} The number of deaths with respect the degree $k$ (social contacts) and  the percentage $c$. The diagram separates the domain into different zones of strategies according to resulting mortality. \textbf{(b)}The expected number of deaths with respect the c\% of . Three representative degrees $k=8, 14, 18$ are shown. \textbf{(c)} Time series of the number of death for constant mixture $a=0.35$ and for degree $k=8, 14, 18$. }

\label{fig:Age3d}            
\end{figure}

\begin{figure}[ht]
\begin{center}
\includegraphics[width=12cm]{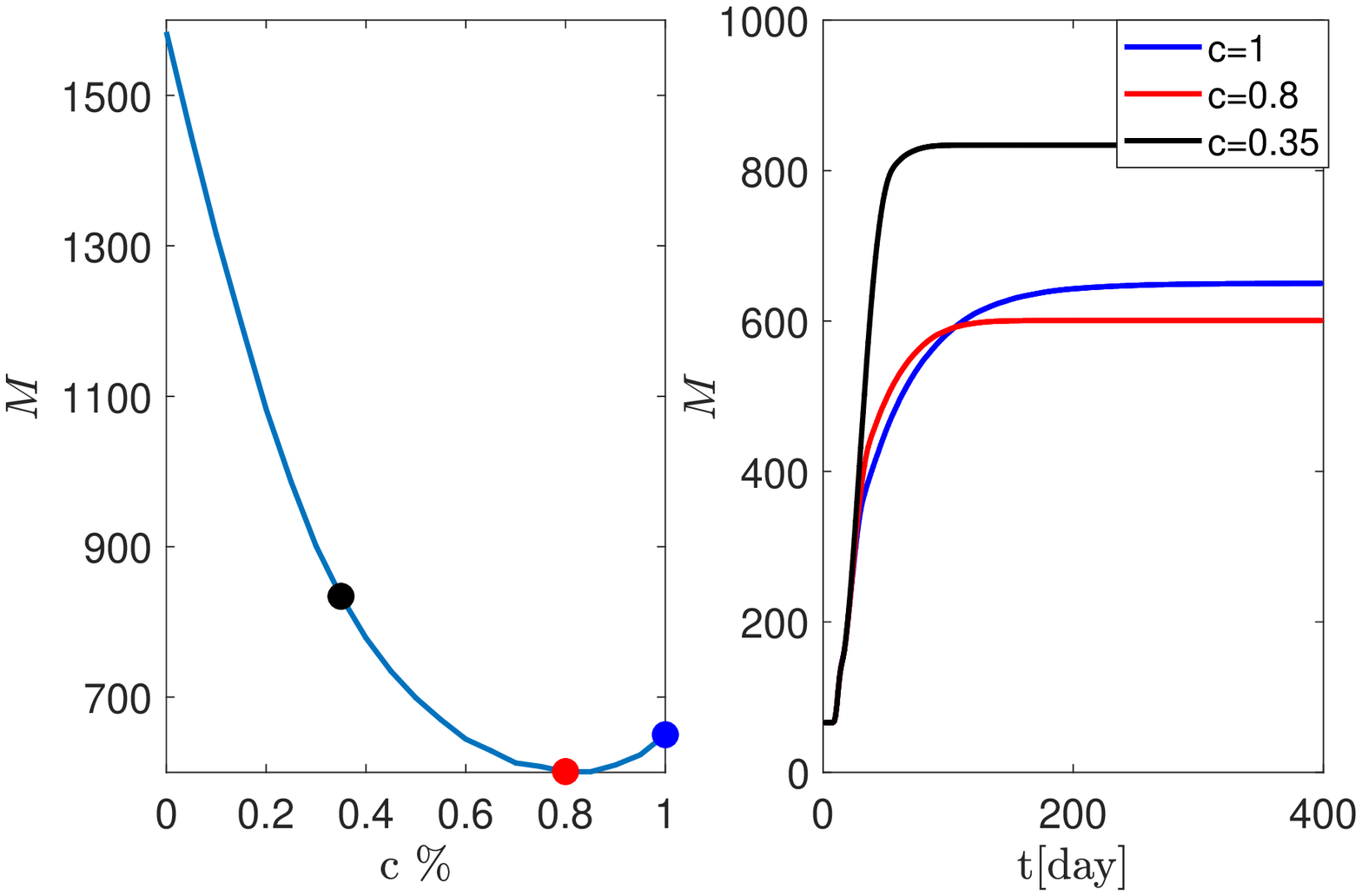}
\end{center}
\caption{Simulation of epidemic model for degree $k=14$. \textbf{(a)} c\% describes the percentage of above 65 year which are vaccinated. Clear there is an infimum at $a=0.8$. instead of $a=1.$ \textbf{(b)} Corresponding curves for $c=0.35$, $c=0.8$ minimum and $c=1$.}
\label{fig:Age3det}            
\end{figure}

\subsection{There exists an optimal age-depened vaccination policy}

The next fundamental question concerns the policies against SARS-CoV-2 pandemic is the age priority during daily vaccinations. The previous simulation results showed in single/double dosing problem that the single dose strategy reduce in optimal way the COVID mortality. Due to this results, we consider the age priority problem for single dose high efficacy vaccination of  $e_0=0.7$, for a vaccination rate 1\% per day. The two groups A and B defined at the begin of this section and the ratio $c\%$ express the priority between the two age groups ($c\%$ gives the proportion of the group A which are daily vaccinated). 

The average dynamics (over the example of 128 realization) with respect the age priority $c\%$ and the connectivity degree are shown in Fig.~\ref{fig:Age3d}(a). This diagram also separates the domain into different zones of strategies according to resulting mortality. For example we have 6 zones (dark blue, blue and light blue, green or cyan and the narrow yellow zones(dark and light yellow) with the two last corresponding to highest mortality.

It is clear that higher values of $c\%$ leads to low mortality, which can be interpreted as high priority vaccination of elder people. It is interesting that can be achieved a low mortality even for high contact rates. For example, in \ref{fig:Age3d}(a) the blue areas of low mortality exist in all range of degrees and importantly, even with degrees over 25 (for $c>0.7$).

The most interesting result is the non monotonic behavior of mortality for the medium connectivity regime with respect to $c$. This is better depicted in \ref{fig:Age3d}(b), for $ k=12,k=18$ degree, where a global minimum is achieved at $c\approx 0.8$. In Fig.~\ref{fig:Age3det} a more detailed illustration is given. The mean degree is held at $k=14$. The expected number of death for  $c=0.4, 0.85$ and 1 are marked with filled circles and resulting  833, 600 and 648 respectively. The best policy for this medium social distancing is at $c=0.85$ and not at $c=1$ which corresponds to "first the older" policy. The importance of best policy assessment, finally, is depicted from the mortality range, where for $k=14$ the death toll  $M_c\in (600, 1582)$.  

\section{Discussion - Conclusion}

In  the context of the COVID-19 pandemic, we studied  the vaccination effect in conjunction with NPI's with a particular focus on social sitancing. A large scale agent-based model SEIRM was used on a complex small-world networks. This combination allowed us to design vaccination strategies (i.e. single or double dosing policy, age stratification) taking into account the impact of social distancing and to prevent the expected number of deaths.   

 A rigorous methodology was developed to estimate the reduction in the infection probability for vaccinated individuals given the population vaccine efficacy. This was achieved by imitating real clinical trials process of efficacy estimation and subsequently by wrapping around the model a numerical solver Newton-Raphson. The proposed methodology has similarities to  Equation Free framework \cite{Proctor14, Kev09} because it  bypasses the need of an equation formulation of efficacy for the Newton-Raphson scheme.

The analysis showed that the optimal vaccination strategy is the single dose deployment, but the expected number of deaths varies significantly under social distancing measures. Deploying the first vaccine dose across the whole population is sufficient to control the epidemic, in terms of COVID-19  deaths, only for mild social distancing measures (low-medium amount of contacts, see Figs\ \ref{fig:Roll500} and \ref{fig:Roll1000}). Remarkably, similar result was obtained from the analysis of \cite{maier21} where using a detailed model of differential equations. They showed that delaying the second vaccine dose is expected to prevent COVID-19 deaths in a four to five digit range. However, their analysis does not consider social distancing measures as the present study.

Using a single dose strategy, the age vaccination priority is investigated. Interestingly, the results have shown the existence of non monotonic behavior for mild social distancing measures. The minimum of the expected number of deaths  resulting from vaccination strategy of 80\% over 65 year old, and the remain percentage from younger groups. The result is robust from single dose vaccine efficacy above 50\%. Similar result obtained in \cite{Matr21} when they used a detailed stratified continuous model and optimization algorithms to estimate optimal vaccine allocation. When minimizing with respect to the number of deaths, they found that for low vaccine efficacy, it is optimal to vaccinate older groups first. In contrast, for higher vaccine efficacy, there is a switch to allocate vaccine to high-transmission individuals (younger people). Additionally, in Buder et al. \cite{Bub21} using again a mathematical model to investigated age-stratified prioritization strategies they reached to similar conclusions, that when  optimization function is the mortality and years of life lost, the minimization - in most of the scenarios - was obtained when the vaccine was prioritized to adults greater than 60 years old \cite{Bub21}. 

Our analysis highlighted the strong dependence of the expected number of deaths on the social distancing interventions. Both Figs.~\ref{fig:Roll500}, \ref{fig:Roll1000} showed the disadvantage of high connectivity (degree$>30$), where even with the optimal strategy of one dose, the expected number of deaths will be almost double compared to the zone of mild connectivity  (i.e. k=15). In addition, mild to higher connectivity is making the mortality of COVID-19 more sensitive to roll out policy. For example in Fig.~\ref{fig:Roll1000} with connectivity ($\approx 20$) the expected deaths ranged from $\approx 600$ ($a\%=0$) to $\approx 1000$ ($a\%=0.7$). The authors in \cite{Moo21} studied multiple scenarios of NPI relaxation and vaccine characteristics. Similar to our results, they showed that vaccination alone is insufficient to contain the outbreak \cite{Moo21}. Specifically, in the absence of NPIs, they reported that even with the optimistic scenario of high vaccination efficacy the reproduction number will be  higher than one ($\approx 1.58$). Furthermore, they reported that after the end of the vaccination program, the expected removal of all NPIs is predicted to lead to more than 20,000 deaths in UK, in the optimal scenario high vaccination efficacy.

The inherent complexity of SARS-CoV-2 pandemic implies that our approach bears limitations with respect to the model assumptions. For example, the model does not distinguishes between asymptomatic and symptomatic infections rate. Furthermore, the simulations are made for the expected number of deaths without considering other possible NPI measures ---(i.e. masks, school regulations, hand washing etc.)---. In addition, side effects which can increase the number of deaths are excluded from model. For example in the case of very restricted NPIs corresponds to very low k degree in our model the following harmful health effects are reported \cite{Ben21}: hunger,drugs overdoses, postponed surgeries and generally missed health services. For example delayed diagnosis of cancer and the suboptimal care may have strong impact for the wider population of cancer patients  \cite{Rich20}. As final concluding remark, we stress that the analysis of the model proposed an alternative to the strong restricted NPI's, since mild social distancing measures (low-medium amount of contacts) simultaneously with optimal vaccine allocation can decrease sufficiently the expected number of deaths.
  
One major future topic will be the use of more realistic structures of network and the impact on the emergent pandemic dynamics. The use of household structure in the connectivity, the dominant heterogeneity factor \cite{Hil19}, will attribute in the realistic modelling of lockdown effect: stronger contacts and hence transmission potential within the family setting \cite{Hil19}. Contributions of detailed demographic data like age structure and age-specific fatality \cite{Dud20} as well as mobility variations during pandemic \cite{BAS21} will also improve our model.

Another research direction is the investigation of the parameter dependence network dynamics (infection rates, efficacy,density of connections). Using equation free, numerical bifurcation tools for multiscale agent based problems \cite{Kev09,Siettos2016,Spil10, Rep10}, will help to predict signs or early warnings for upcoming waves pandemics.

\bibliographystyle{unsrt}  
\bibliography{references.bib} 





\end{document}